\documentstyle[12pt]{article}
\input mssymb
\newtheorem{theorem}{Theorem}[section]
\newtheorem{defn}[theorem]{Definition}

\newtheorem{lemma}[theorem]{Lemma}

\newtheorem{cor}[theorem]{Corollary}
\newtheorem{claim}[theorem]{Claim}

\newtheorem{eple}[theorem]{Example}
\newtheorem{rmk}[theorem]{Remarks}
\newtheorem{dsc}[theorem]{Discussion}
\newtheorem{nota}[theorem]{Notation}
\newtheorem{conv}[theorem]{Convention}

\newsavebox{\indbin}
\savebox{\indbin}{\begin{picture}(0,0)
\newlength{\gnu}
\settowidth{\gnu}{$\smile$}
\setlength{\unitlength}{.5\gnu}
\put(-1,-.65){$\smile$}
\put(-.25,.1){$|$}
\end{picture}}
   
\newcommand{\disp}{\displaystyle}   
\newcommand{\be}{\begin{enumerate}}
\newcommand{\bd}{\begin{defn}}
\newcommand{\bt}{\begin{theorem}}
\newcommand{\bl}{\begin{lemma}}
\newcommand{\ee}{\end{enumerate}}
\newcommand{\ed}{\end{defn}}
\newcommand{\et}{\end{theorem}}
\newcommand{\el}{\end{lemma}}

\newcommand{\la}{\langle}
\newcommand{\ra}{\rangle}
\newcommand{\sub}{\subseteq}

\newcommand{\sm}{\setminus}

\newcommand{\proof}{\noindent{\bf Proof:}\ }
\newcommand{\qed}{\hfill $\Box$}

\newcommand{\Cl}{{\cal L}}

\newcommand{\dom}{\mbox{dom}}

\title{Categoricity over P for first order T or categoricity for $\varphi \in
{\cal L}_{\omega_1\omega}$ can stop at $\aleph_k$ while holding for
$\aleph_{0}, \ldots ,\aleph_{k - 1}$ \thanks{To make Leo happy}}
\author{Saharon Shelah \\ The Hebrew University \\ Department of Mathematics \\
Jerusalem, Israel \thanks{Partially supported by the BSF, grant 323}
\and Bradd Hart \\ University of California at Berkeley \\ 
Department of Mathematics \\ Berkeley, California \thanks{Supported by a grant
from the NSERC}}
\begin{document}
\maketitle
In the 1950's, Los conjectured that if $T$ 
was countable first order 
theory in a language $\cal L$ then if it was categorical in some
uncountable power then it was categorical in all uncountable powers.
In \cite{Morley}, Morley proved this.  Buoyed by this success,  more
general forms of the Los conjecture were considered.

In \cite{ShCat},  Shelah showed that if $T$ was any first order theory
categorical in some power greater than $|T|$ then $T$ was categorical in
all powers greater than $|T|$.  Keisler took up the investigation of
the ${\cal L}_{\omega_1\omega}$ case (see \cite{Keisler}) and gave a
sufficient condition for the Morley analysis to work in this situation.
Unfortunately, this condition was not necessary.  (See the counter-example
due to Marcus, \cite{Marcus})

In \cite{CatLw1w} and \cite{classnon1}, Shelah began the systematic
investigation of the ${\cal L}_{\omega_1\omega}$ case.  In
\cite{classnon1}, he identifies a class of ${\cal L}_{\omega_1\omega}$
sentences which he calls excellent and shows that if an 
${\cal L}_{\omega_1\omega}$ sentence is excellent then the Los conjecture
holds.  (In \cite{thesis}, Hart shows that many other theorems which
are analogs of those for first order theories also hold for excellent
classes.)  Furthermore, he shows that if GCH (or in fact much less)
and $\varphi$ is an ${\cal L}_{\omega_1\omega}$ sentence which is
$\aleph_n$-categorical for all $n \in \omega$ then $\varphi$ is
excellent.

The question which naturally arises is, under suitable set theoretic
assumptions, is categoricity in $\aleph_n$ for $n < k$ sufficient to
prove full categoricity for a sentence in ${\cal L}_{\omega_1\omega}$.

The answer to this question must wait while we introduce another variant
of the Los conjecture.

Suppose ${\cal L}$ is a relational language and $P \in {\cal L}$ is a
unary predicate.  If $M$ is an ${\cal L}$-structure then $P(M)$ is the
${\cal L}$-structure formed as the substructure of $M$ with domain
$\{a : M \models P(a) \}$.  Now suppose $T$ is a complete first order
theory in ${\cal L}$ with infinite models.  Following Hodges, we define

\begin{defn}
T is relatively $\lambda$-categorical if whenever $M$,$N \models T$,
$P(M) = P(N)$, $|P(M)| = \lambda$ then there is an isomorphism
$i : M \rightarrow N$ which is the identity on $P(M)$.

T is relatively categorical if it is relatively $\lambda$-categorical
for every $\lambda$.
\end{defn}

The notion of relative categoricity has been investigated by Gaifman
(\cite{Gaifman}), Hodges (\cite{Hodges} and \cite{relcathod}), Pillay
(\cite{wcat/P}) and Pillay and Shelah (\cite{Stab/P}). In
(\cite{class/p}), Shelah gave a classification under some set theory.

Again the question arises whether the relative $\lambda$-categoricity
of $T$ for some $\lambda > |T|$ implies that $T$ is relatively categorical.

In this paper, we provide an example, for every $k > 0$, of a theory
$T_k$ and an ${\cal L}_{\omega_1\omega}$ sentence $\varphi_k$ so that
$T_k$ is relatively $\aleph_n$-categorical for $n < k$ and $\varphi_k$
is $\aleph_n$-categorical for $n < k$ but $T_k$ is not
relatively $\beth_k$-categorical and $\varphi_k$ is not
$\beth_k$-categorical.

The examples are due to Shelah.  Harrington asked about the
possibility of such examples in Chicago in December, 1985 as he was
not happy with the complexity of the classification.  The examples
provided ${\cal L}_{\omega_1\omega}$ sentences which were categorical
but not excellent and so a proof of this fact was written up in
\cite{thesis}.

The notation used is standard.  $[A]^k$ will stand for all the
$k$-element subsets of the set $A$.  ${\cal P}^-(n)$ is the set of all
subsets of $n$ except $n$ itself.  $\disp{\coprod}$ is used to represent
the direct sum of groups and $\disp{\prod}$ is used to represent the
direct product of groups.  $Z_2$ will represent the two element group.
$2^{<\omega}$ will be used to represent the subgroup of eventually
zero sequences in the abelian group $\displaystyle{\prod_\omega} Z_2$
(written as $2^\omega$).

\section{The Example}

We first describe the example informally.  Fix a natural number $k$
greater than one.  There will be an infinite set $I$ with $K = [I]^k$.
There are constants $c_n$ for $n \in \omega$ and a predicate $R$
containing all of them.  $R$ will be thought of as levels and we will
refer to constants in $R$ as standard levels.  We fix $Z_2$, the abelian
group, $G$, the direct sum of $K$-many copies of $Z_2$ and $H$, the direct
sum of $R$-many copies of $Z_2$.  In addition, all relevant projections
onto $Z_2$ are available to us.  All of this constitutes the $P$-part of
the model.

Outside of this we have two types of objects.  First, for every level
$r \in R$ and every $u \in K$, we have a distinct copy of $G$.  Via some
connection between our fixed copy of $G$ and this one we will be able
to determine the sum of any three elements of $G$ but we will have 
``lost'' the zero.  Second, for every $u \in K$ there will be a
distinct copy of $H$ in which we again have ``lost'' the zero.

We will be interested in the possibility of choosing elements from these
copies of $G$ and $H$ to act as the zero in their respective groups.  We
won't put any more restraints on $G$'s from non-standard levels so
any element will do.   However, for each $n \in \omega$, on the
level corresponding to $c_n$, and for every $u \in K$, there will
be a predicate connecting the copy of $H$ corresponding to u and k
of the copies of $G$ on the $n^{th}$ level.  It will be these predicates
which make or break the categoricity by putting restraints on choices for the
zeroes of the copies of $G$ and $H$.

%

We now wish to fix $k$ for the rest of the paper.

\begin{conv}
$k$ will be a fixed natural number greater than one.
\end{conv}

Now, more formally, we define the language for the example.
\bd \rm $\Cl$ will be the language that consists of
\begin{enumerate}
\item unary predicates $I,K,R,P,G^a,H^a$
\item binary predicates $\in$, $H^b$
\item ternary predicates $\pi$,$\rho$, $+$ and $G^b$
\item a 4-ary predicate $h$
\item a 5-ary predicate $g$
\item a $k+1$-ary predicate $Q_l$ for every $l < \omega$ and
\item constants $c_a$ for every $a \in Z_2 \cup \omega$
\end{enumerate}
\ed

We now describe the standard model on $I$.

\bd \rm If $I$ is an infinite set then the standard model on $I$ denoted
by $M_I$ is the $\Cl$-structure with universe
\[ I \cup [I]^k \cup \omega \cup Z_2 \cup
\coprod_{[I]^k} Z_2 \cup \coprod_\omega Z_2 \cup
\omega \times [I]^k \times \coprod_{[I]^k} Z_2 \cup
[I]^k \times \coprod_\omega Z_2 \]
and where the symbols of $\Cl$ are interpreted as follows:
\begin{enumerate}
\item $I$ is interpreted as $I$, $K$ as $[I]^k$, $R$ as
$\omega$, $G^a$ as $\disp{\coprod_{[I]^k}} Z_2$ and $H^a$ as
$\disp{\coprod_\omega} Z_2$
\item the constants $c_a$ are interpreted as a.  That is, for
example,$R(c_a)$ holds for every $a \in \omega$.
\item $P(x)$ holds iff $x$ is a constant or one of $I(x)$,$K(x)$,$G^a(x)$
or $H^a(x)$ holds.
\item $G^b(l,u,x)$ holds iff $R(l)$, $K(u)$ and $x = (l,u,y)$ for
some $y \in \disp{\coprod_{[I]^k}} Z_2$
\item $H^b(u,x)$ holds iff $K(u)$ and $x = (u,y)$ for some $y \in
\disp{\coprod_\omega} Z_2$
\item $\in (x,y)$ holds iff $I(x)$, $K(y)$ and $x \in y$
\item \label{add1} $+(x,y,z)$ holds iff $x,y$ and $z$ are all in one of $Z_2$, 
$\disp{\coprod_{[I]^k}} Z_2$ or $\disp{\coprod_\omega} Z_2$ and
\mbox{$x + y = z$}.
\item $\pi(u,x,a)$ holds iff $K(u)$, $G^a(x)$ and $x(u) = a$, an
element of $Z_2$.
\item $\rho(l,x,a)$ holds iff $R(l)$, $H^a(x)$ and $x(l) =a$, an
element of $Z_2$.
\item \label{add2} $g(l,u,x,y,z)$ holds iff $R(l)$, $K(u)$, $G^a(x)$,
$y = (l,u,a)$, $z = (l,u,b)$ (so $G^b(l,u,y)$ and $G^b(l,u,z)$) and
\mbox{$b = a + x$}. 
\item \label{add3} $h(u,x,y,z)$ holds iff $K(u)$, $H^a(x)$, $y =
(u,a)$, $z = (u,b)$ (so $H^b(u,y)$ and $H^b(u,z)$) and \mbox{$b = a +
x$}. 
\item \label{add4} $Q_l(x_0, \ldots,x_k)$ holds iff $x_i =
(c_l,u_i,y_i)$  with $G^b(c_l,u_i,x_i)$ for $i < k$ and \mbox{$x_k =
(u_k,z)$} with \mbox{$H^b(u_k,x_k)$} where $u_0, \ldots, u_k$ are 
all the $k$-element subsets of some $(k + 1)$-element subset of $I$ and
\[\sum_{i < k} y_i (u_k) = z(c_l)\]
\end{enumerate}
\ed

{\bf Remarks:} \be \item In the previous definition, all of the direct
sums used in the definition of the universe represent abelian groups.
Hence on the right hand side of items \ref{add1}, \ref{add2},
\ref{add3} and \ref{add4}, the addition mentioned is addition in the
appropriate group.
\item In item \ref{add4}, each $y_i$ is in $\disp{\coprod_{[I]^k}}
Z_2$ and $u_k$ is in $[I]^k$ so $y_i(u_k)$ is in $Z_2$.  $z$ is in
$\disp{\coprod_\omega} Z_2$ and $c_l \in \omega$ so $z(c_l)$ is in
$Z_2$. Hence, the displayed equality is comparing elements of $Z_2$.
\ee

Let's consider some of the sentences in $\Cl$ that the standard model
satisfies.  For a fixed infinite set $I$, $M_I$ satisfies:
\be
\item \label{first} $I$ is an infinite set, $K$ is the collection of
$k$-element subsets of $I$ and $\in$ is the membership relation
between elements of $I$ and elements of $K$.
\item $I,K,R,G^a,H^a$ are disjoint and their union together with
the constants $c_a$ for $a \in Z_2$ form $P$.
\item $R(c_a)$ for every $a \in \omega$.
\item $G^b(l,u,x)$ implies $R(l)$ and $K(u)$ and $H^b(u,x)$ implies
$K(u)$.
\item If $x$ is not in $P$ then either for some $l$ and $u$, $G^b(l,u,x)$ or
for some $u$, $H^b(u,x)$ and for every $l \in R$ and $u,v \in K$, $P$,
$H^b(u,-)$ and $G^b(l,v,-)$ are pairwise disjoint.
\item If $\pi(u,a,z)$ then $K(u)$, $G^a(a)$ and $z$ is one of the constants 
indexed by $Z_2$.
\item If $\rho(l,b,z)$ then $R(l)$, $H^a(b)$ and $z$ is one of the constants 
indexed by $Z_2$.
\item If $g(l,u,a,v,w)$ then $R(l)$, $K(u)$, $G^a(a)$,  $G^b(l,u,v)$ and
$G^b(l,u,w)$.
\item If $h(u,b,x,y)$ then $K(u)$, $H^a(b)$, $H^b(u,x)$ and $H^b(u,y)$.
\item The constants $c_a$ for $a \in Z_2$ together with + have the 
group structure of $Z_2$.
\item + restricted to $G^a$ gives a subgroup of $\disp{\prod_K} Z_2$ which
contains $\disp{\coprod_K} Z_2$ where the projections are given by $\pi$.
\item + restricted to $H^a$ gives a subgroup of $\disp{\prod_R} Z_2$ which
contains $\disp{\coprod_R} Z_2$ where the projections are given by $\rho$.
\item For every $l$ in $R$ and $u$ in $K$, $G^b(l,u,-)$ is non-empty and
for every $l$ in $R$, $u$ in $K$ and $x$ so that $G^b(l,u,x)$ $g(l,u,-,x,-)$
is a bijection from $G^a$ onto $G^b(l,u,-)$.  Moreover,
$g(l,u,x,y,z)$ implies $g(l,u,x,z,y)$ and if $g(l,u,a,x,y)$ and $g(l,u,b,y,z)$
then $g(u,l,a + b,x,z)$ where $a + b$ is the unique $c$ so that
$+(a,b,c)$,.
\item For every $u$ in $K$, $H^b(u,-)$ is non-empty and for every $u$ in $K$
and $x$ so that $H^b(u,x)$, $h(u,-,x,-)$ is a bijection from $H^a$ onto
$H^b(u,-)$.  Moreover, $h(u,x,y,z)$ implies $h(u,x,z,y)$ and if
$h(u,a,x,y)$ and $h(u,b,y,z)$ then $h(u,a + b,x,z)$ where $a + b$ is the
unique $c$ so that $+(a,b,c)$.
\item \label{Q1} If $Q_l(x_0, \ldots,x_k)$ then for $i < k$, for some
$u_i$ in $K$, $G^b(c_l,u_i,x_i)$ and for some $u_k$ in $K$,
$H^b(u_k,x_k)$.  Additionally, $u_0, \ldots, u_k$ are all the
$k$-element subsets of some $(k+1)$-element subset of $I$. If $\sigma$
is a permutation of $k$ then $Q_l(x_{\sigma(0)},\ldots,x_{\sigma(k-1)},x_k)$.
\item \begin{sloppypar} \label{Q2}
If $Q_l(x_0,\ldots,x_k)$, $G^b(c_l,u,x_0)$, $H^b(v,x_k)$,
$G^b(c_l,u,x'_0)$ and $H^b(v,x'_k)$
then $Q_l(x'_0,\ldots,x_k)$ iff the $v$-projection of the unique element
$a$ so that $g(c_l,u,a,x_0,x'_0)$ via $\pi$ is 0 and
$Q_l(x_0,\ldots,x'_k)$ iff the $c_l$-projection of the unique element
$a$ so that $h(v,a,x_k,x'_k)$ via $\rho$ is 0.
\end{sloppypar}
\item \label{Q3} Suppose $l \in \omega$, $u$ is in $K$ and
$i_0,\ldots,i_{n-1}$ are distinct elements of $I$ not in $u$.  For
each $j < n$, let $v^j_i$ for $1 \leq i \leq k$ be a list of the
$k$-element subsets of $u \cup \{ i_j \}$ besides $u$.  If
$G^b(c_l,v^j_i,x^j_i)$ for each $j < n$ and $i < k$ and
$H^b(v^j_k,y_j)$ for every $j < n$ then
\[ \exists x \bigwedge_{j < n} Q_l(x,x^j_1,\ldots,x^j_{k-1},y_j) .\]
\ee

\ref{Q3} actually follows from the previous axioms but it is in the form
that we will use it in section 2. We make the following definition for
the rest of the paper.

\begin{conv}
Let $T$ be the theory in $\Cl$ made up of the sentences enumerated
\ref{first} -- \ref{Q3} above. 
\end{conv}

The standard model satisfies some additional sentences in ${\cal
L}_{\omega_1 \omega}$.  For any infinite set $I$, $M_I$ satisfies:
\be \item $R$ contains only the constants
indexed by $\omega$.
\item $G^a$ is canonically isomorphic to $\disp{\coprod_K} Z_2$.
\item $H^a$ is canonically isomorphic to $\disp{\coprod_\omega} Z_2$.
\ee

\begin{conv}
Let $\varphi$ be the ${\cal L}_{\omega_1 \omega}$ sentence which is
the conjunction of $T$ and the three sentences listed above.
\end{conv}

\begin{sloppypar}
{\bf Remarks:} \be \item $T$ is not complete however we will show that
it is relatively $\aleph_n$-categorical for all $n < k$.
\item $\varphi$ is the Scott sentence of any $M_I$ where $I$ is
countable.  This will follow from section 2.  Note that $\varphi$ has
arbitrarily large models.
\ee
\end{sloppypar}

\section{Categoricity less than $\aleph_k$}

In this section, we show that T is relatively $\aleph_n$-categorical
for all $n < k$.


\begin{defn} 
Suppose $M \models T$, $W \subseteq \omega \times K(M) \cup K(M)$
and $f:W \rightarrow M$.  Then f is called a solution for W if:
\begin{enumerate}
\item $(l,u) \in W$ then $M \models G^b(c_l,u,f(l,u))$
\item if $u \in W$ then $M \models H^b(u,f(u))$ and
\item if $u_0,\ldots,u_k \in K(M)$ are all the $k$-element subsets
of some fixed $(k+1)$-element subset of I(M), $(l,u_i) \in W$ for
all $i < k$ and $u_k \in W$ then
\[ M \models Q_l(f(l,u_0),\ldots,f(u_k)) \]
\end{enumerate}

If $J \subseteq I(M)$ then f is called a J-solution if it is a
solution for $\omega \times [J]^k \cup [J]^k$.  f is called a
solution if it is an I(M)-solution
\end{defn}

{\bf Remark:} Note that the standard model for any $I$ has a solution.
Hence $T$ (and $\varphi$) has arbitrarily large models with solutions.

\bl \label{solution}
If $M,N \models T$, both M and N have solutions and $P(M) = P(N)$
then $M \cong N$ over $P(M)$.
\el

\proof Suppose $f_M$ is a solution for $M$ and $f_N$ is a solution for
$N$.  We are really interested in those $G^b(u,M)$ and $G^b(l,u,N)$
where $l$ is one of the constants in $R$.  However, we must
accommodate all $l$ in $R$.  Let \[R^* = R(M) \sm \{ c_l : l \in \omega
\}.\]  Extend $f_M$ and $f_N$ to include $R^* \times K(M)$ (= $R^*
\times K(N)$) in their domains so that
\[ M \models G^b(l,u,f_M(l,u)) \mbox{ and }  N \models
G^b(l,u,f_N(l,u)) \] for all $(l,u) \in R^* \times K(M)$.
Let $j$ be a partial function from $M$ to $N$ so that $j$
restricted to $P(M)$ is the identity, for every $u$, $j(f_M(u)) =
f_N(u)$ and for every $l$ and $u$, $j(f_M(l,u)) = f_N(l,u)$.  We want
to extend $j$ to a function from $M$ to $N$.

If $x \in M$ so that $M \models G^b(c_l,u,x)$ then there is a unique
$a$ so that 
\[M \models g(c_l,u,a,f_M(l,u),x).\]
There is a unique $y \in N$
so that 
\[N \models g(c_l,u,a,f_N(l,u),y).\]
Extend $j$ so that $j(x) = y$.

We do a similar thing when $x \in M$, $M \models G^b(l,u,x)$ and $l
\in R^*$.

If $x \in M$ so that $M \models H^b(u,x)$ then there is a unique $a$
so that \[M \models h(u,a,f_M(u),x).\] There is a unique $y \in N$ so
that \[N \models g(u,a,f_N(u),y).\] Extend $j$ so that $j(x) = y$.

Using the fact that $M$ and $N$ satisfy $T$, it is not hard to show
that $j$ defines a function from $M$ onto $N$.  We want to show that
it is an isomorphism.  We'll check the hardest predicate, $Q_l$.

Suppose $M \models Q_l(x_0,\ldots,x_k)$ where \[M \models
G^b(c_l,u_i,x_i) \mbox{ for $i < k$ and } M \models H^b(u_k,x_k).\]  Choose
$a_i$ for $i \leq k$ so that $M \models g(c_l,u_i,a_i,f_M(l,u_i),x_i)$
for $i < k$ and \[M \models h(u_k,a_k,f_M(u_k),x_k).\]  We know \[M
\models Q_l(f_M(l,u_0),\ldots,f_M(u_k))\] since $f_M$ is a solution. 
Suppose \[M \models \pi(u_i,a_i,z_i) \mbox{ for $i < k$ and } M \models
\rho(c_l,a_k,z_k).\]  Then by using axioms \ref{Q1} and \ref{Q2} of
$T$, we conclude that
\[\sum_{i<k} z_i = z_k \] where the sum takes place in $Z_2$ and we
identify the constants indexed by $Z_2$ with the elements they
represent.

Since $P(M) = P(N)$, this happens in $N$ as well and since $N \models
T$, we unravel the fact that $f_N$ is a solution so $N \models
Q_l(f_N(l,u_0),\ldots,f_N(u_k))$ to conclude that $N \models
Q_l(y_0,\ldots,y_k)$ where $y_i = j(x_i)$ for $i \leq k$.

A completely symmetric argument shows that if $N \models
Q_l(j(x_0),\ldots,j(x_k))$ then $M \models Q_l(x_0,\ldots,x_k)$ so $j$
is an isomorphism. \qed

\bl \label{iso} Suppose $M \models T$.
\be \item If $M$ is countable then $M$ has a solution.
\item If $A \subseteq B \subseteq I(M)$, $B$ is
countable and $f$ is an $A$-solution then $f$ can be extended to a
$B$-solution. \ee
\el

\proof The first follows from second so we will prove the second. 

Choose $f'$ so that $f \sub f'$ and $\dom(f') = \dom(f) \cup [B]^k$
where if $u \not \in [A]^k$ then $M \models H^b(u,f'(u))$ and
otherwise $f'(u)$ is arbitrary.

$f'$ is a solution on its domain.  To see this, note that if
$i_0,\ldots,i_k \in B$ and $i_0 \not \in A$ then since $k > 1$, at
least two $k$-element subsets of $\{i_0,\ldots,i_k\}$ are not in
$[A]^k$.  Hence, $f'$ is a solution on its domain vacuously.

Now enumerate $\omega \times ([B]^k \sm [A]^k)$ as $\{ \la l_i, u_i
\ra : i \in \omega \}$.  We will define an increasing chain of
functions $f_n$ so that \be \item $f_0 = f'$,
\item $\dom(f_n) = \dom(f') \cup \{ \la l_i,u_i\ra : i < n \}$ and
\item $f_n$ is a solution on its domain. \ee
If we accomplish this then $\bigcup f_n$ will provide a $B$-solution
extending $f$.

Suppose we have defined $f_n$.  We need to choose an $a$ so that $M 
\models G^b(c_{l_n},u_n,a)$ and which will be compatible with the demands
of being a solution.

Say that a $(k+1)$-element subset $v$ of $B$ puts a constraint on
$u_n$ if $u_n \sub v$ and $k-1$ of the $k$-element subsets of $v$, say
$w_1,\ldots,w_{k-1}$, are such that $\la l_n, w_i \ra \in \dom(f_n)$
for $i < k$.  Note that since $u_n \not \sub A$, at least one of these
$w_i$'s must also not be a subset of $A$.

Now since only finitely many elements are enumerated before $\la l_n,
u_n \ra$, there are only finitely many $(k+1)$-element subsets of $B$
which put a constraint on $u_n$.  This is exactly the situation that
axiom \ref{Q3} of $T$ was designed for so we can find an $a$ so that
$f_{n+1} = f_n \cup \{ \la \la l_n, u_n \ra, a \ra \}$ is a solution
on its domain. \qed

\begin{cor} 
$\varphi$ is a complete $\Cl_{\omega_1\omega}$ sentence.
\end{cor}

\proof To see this, it suffices to see that if $M$ and $N$ are
countable models of $\varphi$ then $M \cong N$.  But since $M$ and $N$
are models of $\varphi$, $P(M)$ and $P(N)$ are uniquely determined by
$\varphi$ so we may assume that $P(M) = P(N)$.  By lemma \ref{iso},
$M$ and $N$ have solutions and hence by lemma \ref{solution}, $M \cong N$.
\qed

\begin{defn} Suppose $M \models T$, $A_\emptyset \sub I(M)$ and $a_0,\ldots,a_{m-1}$ are distinct elements of $I(M) \sm A_\emptyset$.
$\la A_s , f_s : s \in {\cal P}^-(m) \ra$ is a compatible $\aleph_n -
{\cal P}^-(m)$-system of solutions if
\begin{enumerate}
\item $\bigcup_{s \in {\cal P}^-(m)} A_s = A_\emptyset \cup \{a_0,\ldots,
a_{m-1}\}$, $|A_\emptyset| \leq \aleph_n$ and $A_s = A_\emptyset
\cup \{a_t : t \in s \}$ for every $s \in {\cal P}^-(m)$.
\item $f_s$ is a $A_s$-solution for every
$s \in {\cal P}^-(m)$
\item for every $s,t \in {\cal P}^-(m)$ 
if $s \subseteq t$ then $f_s \subseteq f_t$
\end{enumerate}
\end{defn}

Using the notation from the definition, suppose $\la  A_s , f_s : s
\in {\cal P}^-(m) \ra$ is a compatible $\aleph_0-{\cal P}^-(m)$-system
with $m < k$.  If \[u \in [\bigcup_{s \in {\cal P}^-(m)} A_s]^k
\setminus \bigcup_{s \in {\cal P}^-(m)}[A_s]^k\] then $\{a_0, \ldots,a_{m-1}\}
\sub u$.  Since $m < k$, there is $b \in u \setminus \{a_0, \ldots,a_{m-1}\}
\sub u$.  If $c \in \bigcup_{s \in {\cal P}^-(m)} A_s \setminus u$ then \[(u
\setminus \{b\}) \cup \{c\} \not \in \bigcup_{s \in {\cal P}^-(m)}[A_s]^k.\]
Hence, if $u \sub v$ where $v$ is any $(k+1)$-element subset of
$\bigcup_{{\cal P}^-(m)} A_s$ then there is a $k$-element subset $u'
\sub v$, $u \neq u'$ so that $u' \not \in \bigcup_{{\cal P}^-(m)}
A_s$ as well.  Using this observation and a proof similar to the proof of
lemma \ref{iso}, we obtain

\bl
\label{ctblep-m}
If $\la  A_s , f_s : s \in {\cal P}^-(m) \ra$ is a compatible $\aleph_{0} -
{\cal P}^-(m)$-system with $m < k$ then there is $\bigcup_{s \in
{\cal P}^-(m)} A_s$-solution f so that $f_s \subseteq f$ for every $s
\in {\cal P}^-(m)$
\el

We use this as the base step in the following lemma

\bl
\label{p-m}
If $\la  A_s , f_s : s \in {\cal P}^-(m) \ra$ is a compatible $\aleph_n - {\cal P}^-(m)
$-system with $m + n < k$ then there is $\bigcup_{s \in {\cal P}^-(m)} A_s$-solution f
so that $f_s \subseteq f$ for every $s \in {\cal P}^-(m)$
\el

\proof  We prove this by induction on $n$.  If $n = 0$ then this is just
lemma \ref{ctblep-m}.  Suppose $n > 0$ and $A_s = A_\emptyset \cup \{b_t :
t \in s \}$.  Enumerate $A_\emptyset$,
$\la a_\beta : \beta < \aleph_n\ra$ and let $A_\emptyset^\alpha = \{a_\beta :
\beta < \alpha\}$.   Now define $A_s^\alpha = A_\emptyset^\alpha \cup
\{ b_t : t \in s \}$ for every $s \in {\cal P}^-(m)$ and let $f_s^\alpha$ be
the restriction of $f_s$ to an $A_s^\alpha$-solution.
We wish to define $g_\alpha$ for every $\alpha < \aleph_n$ so that
\begin{enumerate}
\item $g_\alpha$ is a  $\bigcup_{s \in {\cal P}^-(m)} A_s^\alpha$-solution extending
$f_s^\alpha$ for every $s \in {\cal P}^-(m)$
\item $g_\alpha \subseteq g_\beta$ for $\alpha < \beta < \aleph_n$
\end{enumerate}

Clearly, if we accomplish this then $\bigcup_{\alpha < \aleph_n} g_\alpha$
is the sought after solution.  But by taking unions at limit ordinals and
using the induction hypothesis at successors we can easily satisfy these
two conditions so we are done. \qed

\bl
\label{extend}
If $M \models T$ and $A \subseteq B \subseteq I(M)$ with $|B| < \aleph_{k-1}$
and $f$ is an A-solution then $f$ can be extended to a B-solution.
\el

\proof Without loss of generality, $B = A \cup \{b\}$
We prove this lemma by induction on the cardinality of $A$.  If
$A$ is countable then this is just lemma \ref{iso}.  If $|A| = \aleph_n$
with $n > 0$ then enumerate $A$ as $\la a_\beta : \beta < \aleph_n\ra$
and let $A_\alpha = \{a_\beta : \beta < \alpha \}$
Let $f_\alpha$ be the 
restriction of $f$ to an $A_\alpha$-solution.
By induction, we define  $A_\alpha \cup \{b\}$-solutions $g_\alpha$ extending
$f_\alpha$.  If we have defined  $g_\alpha$,
we use lemma \ref{p-m} in the case $m = 2$ to extend $g_\alpha \cup f_{\alpha
+1}$ to a $A_{\alpha + 1} \cup \{b\}$-solution.  At limits we take unions
and $\bigcup_{\alpha < \aleph_n} g_\alpha$ is a $B$-solution extending
f. \qed

\bt \label{cattheo}
If $M \models T$ and $|M| < \aleph_k$ then M has a solution.
\et

\proof By induction on the cardinality of $M$.  If $M$ is
countable then this is lemma \ref{iso}.
If $|M| = \aleph_n$ with $n > 0$ then we can choose $N$,
$N \prec M$ with $|N| < \aleph_n$.  By induction, $N$ has a
solution and by using lemma \ref{extend} repeatedly, we can
extend it to a solution for $M$. \qed

\begin{cor} \be \item $T$ is relatively $\aleph_n$-categorical for all
$n < k$.
\item  $\varphi$ is $\aleph_n$-categorical for all $n < k$. 
\ee
\end{cor}

\proof 1. Suppose $M$ and $N$ are models of $T$, $P(M) = P(N)$ and
$|P(M)| = \aleph_n$ for some $n < k$. It follows that $|M| = |N| =
\aleph_n$.  By theorem \ref{cattheo}, $M$ and $N$ have solutions and
so by lemma \ref{solution}, $M \cong N$.

2. Suppose $M$ and $N$ are models of $\varphi$ and $|M| = |N| =
\aleph_n$ for some $n < k$.  $P(M)$ is uniquely determined by $I(M)$
and $P(N)$ is determined by $I(N)$.  $|M| = |I(M)|$ so we may assume
that $P(M) = P(N)$ and it follows then that $M \cong N$ by theorem
\ref{cattheo} and lemma \ref{solution}. \qed
\section{The Failure of Full Categoricity}

In this section, we show that $\varphi$ is not fully categorical.

Suppose $M \models \varphi$ and $I = I(M)$.  Without loss of
generality, we may assume that $K(M) = [I]^k$, $R(M) = \omega$,
$G^a(M) = \disp{\coprod_K} Z_2$ and $H^a = \disp{\coprod_\omega} Z_2$.
Further, we may assume that the constants $c_l = l$ for  $l \in
\omega$ and $c_a = a$ for $a \in Z_2$.  $\pi, \rho$ and $+$ can also
be assumed to be as in the standard model $M_I$.

\bl \label{solgodown} If $M,N \models \varphi$, $M \sub N$ and $N$ has
a solution then $M$ has a solution.
\el

\proof Suppose that f is a solution for $N$.  Fix some $g: 
\omega \times K(M)
\rightarrow M$ so that \[ M \models  G^b(l,u,g(l,u)) \mbox{ for every
} l \in \omega \mbox{ and } u \in K(M).\]  For $u \in K(M)$,
let $c_{l,u}$ be such that \[N \models g(l,u,c_{l,u},g(l,u),f(l,u)).\]
Choose $d_{l,u}$ so that for every $v \in K(M)$ and $y \in Z_2$
\[ M \models \pi(v,d_{l,u},y) \mbox{ iff } M \models \pi(v,c_{l,u},y)
.\] Define $f':\omega \times K(M) \cup K(M) \rightarrow M$ so that
$f'(u) = f(u)$ for every $u \in K(M)$ and if $l \in \omega$ and $u \in
K(M)$ then $f'(l,u) = z$ where $M \models g(l,u,d_{l,u},g(l,u),z)$. To
check that $f'$ is a solution for $M$, suppose $v$ is a $k+1$-element
subset of $I(M)$ and $u_0,\ldots,u_k$ are all the $k$-element subsets
of $v$.  Fix $l \in \omega$.  \[ N \models
Q_l(f(l,u_0),\ldots,f(u_k)).\] From above, we have \[ N \models
g(l,u_i,c_{l,u_i},d_{l,u_i},f(l,u_i),f'(l,u_i)) \mbox{ for } i < k\]
and by the choice of $d_{l,u}$, \[(c_{l,u_i} + d_{l,u_i})(u_k) = 0
\mbox{ for all } i < k \] hence $M \models
Q_l(f(l,u_0),\ldots,f(u_k))$. \qed

\bl \label{extension} If $M \models \varphi$ and $\kappa > |M|$ then
there is $N \models \varphi$ so that $|N| = \kappa$ and $M \sub N$.
\el 

\begin{sloppypar}
\proof Let $I(N)$ be the disjoint union of $I(M)$ and $\kappa$.  From
our discussion at the beginning of the section, this defines the
$P$-part of $N$.  $P(M)$ will be subset of $P(N)$ except for $G^a(M)$.
The small technical point here is that we have identified $G^a(N)$
with $\disp{\coprod_{K(N)}} Z_2$.  We will identify $x \in G^a(M)$ with
$x' \in G^a(N)$ where $x'(u) =x(u)$ for all $u \in K(M)$ and $x'(u) =
0$ for all \mbox{$u \in K(N) \sm K(M)$.}  In this way, we embed $P(M)$ into
$P(N)$.
\end{sloppypar}

Let's consider the other predicates.  If $u \in K(M)$ then let
$H^b(u,N) = H^b(u,M)$.  If $u \in K(N) \sm K(M)$, let $H^b(u,N) = 2^{<
\omega}$.  It is clear how to define $h$ for $N$ in a fashion
appropriate for $\varphi$.

Let $J = \disp{\coprod_{K(N) \sm K(M)}} Z_2$.  If $u \in K(M)$ and $l
\in \omega$ then let $G^b(l,u,N) = G^b(l,u,M) \times J$ and identify $x
\in G^b(l,u,M)$ with $(x,0)$ where $0$ is the identity in $J$.  If $u
\in K(N) \sm K(M)$, let $G^b(l,u,N) = \disp{\coprod_{K(N)}} Z_2$. We
leave it to the reader to define a reasonable $g$.

It remains to define $Q_l$ on $N$ for each $l \in \omega$.  Fix an
arbitrary function $f:K(M) \rightarrow M$ so that
\[M \models H^b(u,f(u)) \mbox{ for all } u \in K(M).\]
$f$ is needed only in case 3 below.  Suppose $v$ is a $k+1$-element
subset of $I(N)$ and $u_0,\ldots,u_k$ are all the $k$-element subsets
of $v$.  Note that either $v \sub I(M)$ or at most one of the $u_i$'s
is a subset of $I(M)$.  Further suppose $x_i \in G^b(l,u_i,N)$ for $i
< k$ and $x_k \in H^b(u_k,N)$.  There are a number of cases:
\be \item $u_i \in K(M)$ for all $i$.  Then $x_i = (x_i',a_i)$ for
some $x_i' \in G^b(l,u_i,M)$ and $a_i \in J$ for $i < k$.  Since $u_k
\in K(M)$, let 
\[Q_l(x_0,\ldots,x_k) \mbox{ hold in N iff } M \models
Q_l(x_0',\ldots,x_{k-1}',x_k).\] 
\item For only one $j < k$, $u_j \in K(M)$.  $x_j = (x_j',a_j)$ for
some $a_j \in J$.  Let
\[Q_l(x_0,\ldots,x_k) \mbox{ hold in N iff } 
\sum_{i < k} x_i(u_k) = x_k(l) \]
where $x_j(u_k)$ means $a_j(u_k)$.
\item Only $u_k \in K(M)$.  Choose $c$ so that $M \models
h(u_k,c,x_k,f(u_k))$.  Let
\[Q_l(x_0,\ldots,x_k) \mbox{ hold in N iff } M \models
\sum_{i < k} x_i(u_k) = c(l).\]
\item If none of the $u_i$'s are in $K(M)$ then
\[Q_l(x_0,\ldots,x_k) \mbox{ hold in N iff } 
\sum_{i < k} x_i(u_k) = x_k(l). \]
\ee

It is not hard to see that $N$ defined in this way is a model of
$\varphi$ and with the appropriate identifications, $M \sub N$. \qed

\begin{cor} \label{cor1} If $\varphi$ is not $\lambda$-categorical
then it is not $\kappa$-categorical for any $\kappa > \lambda$.
\end{cor}

\proof Any two models of $\varphi$ of cardinality $\lambda$ have
isomorphic $P$-parts.  Hence if $\varphi$ is not $\lambda$-categorical
there must be $M \models \varphi$, $|M| = \lambda$ so that $M$ does
not have a solution.

By lemma \ref{extension}, we can find $N \models \varphi$ and $M \sub
N$ so that $|N| = \kappa$.  If $\varphi$ is $\kappa$-categorical then
$N$ has a solution since there is a model of $\varphi$ of cardinality
$\kappa$ with a solution.  But then by lemma \ref{solgodown}, $M$ has
a solution which is a contradiction.  Hence $\varphi$ is not
$\kappa$-categorical. \qed

\begin{defn}
Suppose $M \models \varphi$ and $i_0,\ldots,i_k$ are distinct elements
of I(M).  Let $A = \omega \times ([\{i_0,\ldots,i_k\}]^k \setminus
\{i_1,\ldots,i_k\})$ and $f$ be a function with domain containing $A$ so that
\[M \models G^b(l,u,f(l,u)) \mbox{ for all } (l,u) \in A.\]  Let
\[x^j_l = f(l,\{i_0,\ldots,i_{j-1},i_{j+1},\ldots,i_k\}) \mbox{ for $j
\neq 0$ and $l < \omega$}\] and choose $ y \in H^b(\{i_1,\ldots,i_k\},M)$.
Define a function g as follows:
\[ g(l) = \left\{ \begin{array}{ll}
		0 & \mbox{if $M \models Q_l(x^0_l,\ldots,x^{k-1}_l,y)$} \\
		1 & \mbox{otherwise}
	\end{array}
	\right. \]
The invariant for $i_0,\ldots,i_k$ via f is $g + 2^{< \omega}$,
a coset of $2^{< \omega}$ in the abelian group $2^\omega$.
\end{defn}

\bl \label{0inv} The definition of invariant given above is independent of the
choice of $y$.
\el

\proof Use the notation of the definition.  Choose any $y'$ so that
\[M \models H^b(\{i_1,\ldots,i_k\},y').\]
Let $c \in H^a(M)$ be such that \[M \models
h(\{i_1,\ldots,i_k\},c,y,y').\] Let
\[ g'(l) = \left\{ \begin{array}{ll}
		0 & \mbox{if $M \models Q_l(x^0_l,\ldots,x^{k-1}_l,y')$} \\
		1 & \mbox{otherwise}
	\end{array}
	\right. \]
Now $g'(l) = g(l) + c(l)$ for all $l \in \omega$ and $c \in 2^{<
\omega}$ so $g' + 2^{< \omega} = g + 2^{< \omega}$. \qed

If $m \in \omega$ and $f,g$ are functions with the same domain
define the relation $\sim_m$ by \[f \sim_m g \mbox{ iff } |\{x : f(x)
\neq g(x) \}| < \aleph_m.\] 

\bd Suppose $M \models \varphi$, $I \subseteq I(M)$ and 
$i_1,\ldots,i_k$ are distinct elements of  $I(M) \setminus I$.
Let f be a function with domain that contains \[\omega \times ([I \cup
\{i_1,\ldots,i_k\}]^k \setminus \{i_1,\ldots,i_k\}) \]
so that \[ M \models G^b(l,u,f(l,u)) \mbox{ for all } (l,u) \in A.\]
The $0$-invariant for $I,i_1,\ldots,i_k$ via $f$
is the function $g$ with domain $I$ so that $g(a) =$ the invariant for
$a,i_1,\ldots,i_k$ via f.

Suppose $0 < m < k$, $I \subseteq I(M)$ and 
$i_1,\ldots,i_{k-m}$ are distinct elements of \mbox{$I(M) \setminus I$}
and f is a function whose domain contains \[ A = \omega \times ([I \cup
\{i_1,\ldots,i_{k-m}\}]^k \setminus 
\{u : \{i_1,\ldots,i_{k-m}\} \subseteq u\}) \] so that
\[ M \models G^b(l,u,f(l,u)) \mbox{ for all } (l,u) \in A.\]
Let $I_0 \subseteq \ldots \subseteq
I_{m-1} \subseteq I$ where $|I_i| = \aleph_i$.  Choose a function $f'$
so that the domain of $f'$ contains
\[B = \omega \times ([I_{m-1} \cup \{i_1,\ldots,i_{k-m}\}]^k) ,\]
$f'(l,u) \in G^b(l,u,M)$ and $f'$ and f agree on their common domain.

The m-invariant for $I,i_1,\ldots,i_{k-m}$  via
$I_0,\ldots,I_{m-1}$ and f is the $\sim_m$-class of the function
h with domain $I \setminus I_{m-1}$ where
h(a) = the ($m-1$)-invariant for $I_{m-1}$ and $a,i_1,\ldots,i_{k-m}$
via $I_0,\ldots,I_{m-2}$ and $f' \cup f$.
\ed

\bl The definition of $m$-invariant above is independent of the choice
of $f'$.  \el

\proof Note that by lemma \ref{0inv}, the definition of $0$-invariant
is well-defined.  Use the notation of the definition for $m$-invariant
for $m > 0$.  Choose any other applicable $f''$.  Let
\[ \begin{array}{ccl}
C & = & \bigcup \{ v : \exists u \in K(M), l < \omega, c \in G^a(M) 
	\mbox{so that } (l,u) \in B,\\
& & \;\;\;\; c(v) \neq 0 \mbox{ and} M \models g(l,u,c,f'(l,u),f''(l,u)) \}.
\end{array} \]
$|C| \leq \aleph_{m-1}$ since $|B| = \aleph_{m-1}$ and if $a \in I \sm
(I_{m-1} \cup C)$ then the value of $h(a)$ is not effected by the
choice of $f''$ instead of $f'$.  Hence the $\sim_m$-class of $h$ is
well-defined. \qed

Suppose that $I$ is an infinite set and $g : [I]^k \rightarrow
2^\omega / 2^{< \omega}$.  We will define the canonical structure
$M_g$ on $I$ via $g$. 

The $P$-part of $M_g$ is the same as $M_I$.  Moreover, so are the
predicates $G^b$ and $g$.  However, $H^b(u,M_g) = \{ u \} \times g(u)$
for all $u \in [I]^k$.  We modify $h$ so that
\[h(u,x,(u,y),(u,z)) \mbox{ holds in } M_g \mbox{ iff } x + y = z \]
where the addition takes place in $2^\omega$. (Note $2^{< \omega} \sub
2^\omega$.)

The definition of $Q_l$ is identical to the one for $M_I$.  It is not
hard to show that $M_g$ satisfies $\varphi$.

\bt \label{main} Let $\lambda$ be the least cardinal such
that $\lambda^{\aleph_{k-1}} < 2^\lambda$. $\varphi$ is not
categorical in $\lambda$.  In fact, there are $2^\lambda$ many
non-isomorphic models of $\varphi$ of cardinality $\lambda$.

{\em {\bf Remark:} Note that $\aleph_{k-1} < \lambda \leq 2^{\aleph_k}$.}
\et

\proof  Let $B_0 = \{ f_a : a \in 2^\omega / 2^{< \omega}\}$ where
$f_a: \aleph_0 \rightarrow  2^\omega / 2^{< \omega}$ so that $f_a(i) =
a$ for all $i \in \aleph_0$.  Define $B_m$ inductively for $0 < m <
k-1$.   Suppose we have defined
$B_{m-1}$. Let $C = \{ h : h : \aleph_m \setminus \aleph_{m-1}
\rightarrow B_{m-1} \}$. 
Let $B_m$ be a maximal collection of $\not \sim_m$-equivalent
elements in $C$.  It is not hard to show that $|B_m| = 2^{\aleph_m}$.

Fix $A \subseteq B_{k-2}^{\aleph_{k-1}}$ of size $\lambda$.  We
wish to define a structure $M^A$ in such a way as to be able to recover
$A$.  Let $I_A = \aleph_{k-2} \cup \aleph_{k-1} \times \aleph_{k-1}
\cup A$.  Choose $g_A : [I_A]^k \rightarrow  2^\omega / 2^{< \omega}$
so that if  $i_m \in \aleph_m \setminus
\aleph_{m-1}$ for $0 < m < k-1$, $\alpha, \beta < \aleph_{k-1}$
and $a \in A$ then $g(\{a,(\alpha,\beta),i_{k-2},\ldots,i_1\}) = 
a(\alpha)(i_{k-2})\ldots(i_1)$ and otherwise $g(u)$ is arbitrary.
Let $M^A$ be the canonical structure on $I_A$ via $g_A$.

We try to recover $A$ by looking at $(k-2)$-invariants.  We need to
fix certain functions for the rest of the argument.  Let \[\bar f:
\omega \times K(M^A) \rightarrow M^A\] be defined so that $\bar f(l,u)
= (l,u,0)$ where $0$ is the identity element in $\disp{\coprod_{K(M^A)}
Z_2}$.  Remember that $(l,u,0)$ is a member of $G^b(l,u,M^A)$.  Let
$f$ be the restriction of $\bar f$ to $\omega \times [\aleph_{k-2}
\cup \aleph_{k-1} \times \aleph_{k-1}]^k$ and let $h$ be the
restriction of $\bar f$ to $\omega \times [\aleph_{k-2} \cup A]^k$.

\begin{claim} Suppose $m < k-1$ and $i_j \in \aleph_j \sm
\aleph_{j-1}$ for $m < j < k-1$.  The $m$-invariant for
$\aleph_m,i_m,\ldots,i_{k-2},(\alpha,\beta),a$ via
$\aleph_0,\ldots,\aleph_{m-1}$ and $\bar f$ is the $\sim_m$-class of
$a(\alpha)(i_{k-2})\cdots(i_{m+1})$. (If $m=0$ then
$a(\alpha)(i_{k-2})\cdots(i_1)$ is the $0$-invariant.)
\end{claim}

\proof Notice that $\bar f$ contains all possible domains required for
calculating invariants.  $\bar f$ essentially chooses the zero in all
the $G^b(l,u,M^A)$'s.

We prove this claim by induction on $m$.  Suppose the notation is as
it is in the claim.  Choose \[y \in a(\alpha)(i_{k-2})\cdots(i_1) =
H^b(u,M^A)\] where $u = \{i_1,\ldots,i_{k-2},(\alpha,\beta),a\}$.

Since $\bar f$ chooses the zero in all $G^b(l,u,M^A)$'s, the value
$y(l)$ determines the truth value of the appropriate instance of
$Q_l$.  This is independent of the choice of $i_0 \in \aleph_0$ so the
$0$-invariant is $a(\alpha)(i_{k-2})\cdots(i_1)$.

The induction step is similar. \qed

A consequence of the claim is that if $a \in A$ and $\alpha,\beta <
\aleph_{k-1}$ then the ($k-2$)-invariant for
$\aleph_{k-2},(\alpha,\beta),a$ via 
$\aleph_0,\ldots,\aleph_{k-3}$ (if $k > 3$) and $f \cup h$ is the
$\sim_{k-2}$-class of $a(\alpha)$.  The domain of $h$ is too large
however to allow us to say we have captured $a$.

So suppose we use some $h'$ instead of $h$ which agrees with $f$ on their
common domain.  Then for any $a \in A$, the value of at most
$\aleph_{k-2}$ many of the ($k-2$)-invariants calculated above
would be effected.  Hence to recover $a(\alpha)$, for every $\beta <
\aleph_{k-1}$, calculate the ($k-2$)-invariant for
$I_{k-2},(\alpha,\beta),a$ via $\aleph_0,\ldots,\aleph_{k-3}$ and $f
\cup h'$ for any $h'$.  All but at most $\aleph_{k-2}$ of the
$(k-2)$-invariants will agree and this $(k-2)$-invariant will be the
$\sim_{k-2}$-class of $a(\alpha)$.

So by fixing $\aleph_{k-2} \cup \aleph_{k-1} \times \aleph_{k-1}$
and $f$ we are able to recover $A$.  We have fixed $\aleph_{k-1}$ elements
then and there are $2^\lambda$ many possible $A$'s, so $2^\lambda$
many of the $M^A$'s are non-isomorphic since $\lambda^{\aleph_{k-1}}
< 2^{\lambda}$. \qed

\begin{cor} \be \item $\varphi$ is not $2^{\aleph_{k-1}}$-categorical.
\item $T$ is not relatively categorical. \ee
\end{cor}

\proof The first is obvious from theorem \ref{main}, the remark after
it and corollary \ref{cor1}.  To see the second, notice that all the
models built in the proof of theorem \ref{main} have isomorphic
$P$-parts and are models of $T$.  Hence $T$ is not relatively
categorical. \qed

\end{document}